# Par-delà le théorème de cocyclicité de Conway : généralisation et alternatives

Par David Pouvreau[1]

**Résumé**
Le fameux théorème de cocyclicité de John H. Conway est ici reconsidéré au moyen d'un paramétrage de la configuration triangulaire associée avec des triplets quelconques $(\alpha; \beta; \gamma)$ de réels. Ce théorème, qui correspond alors au cas $(\alpha; \beta; \gamma) = (1; 1; 1)$, est généralisé en établissant qu'il existe une famille infinie de tels triplets tels que sa conclusion demeure. La configuration « anti-Conway » correspondant au cas $(\alpha; \beta; \gamma) = (-1; -1; -1)$ est aussi étudiée : le théorème de concourance de droites de Xavier Dussau est redémontré et complété par un autre théorème de cocyclicité. Il est aussi démontré qu'il existe en général un unique triplet $(\alpha; \beta; \gamma) \neq (-1; -1; -1)$ fonction des longueurs des côtés du triangle considéré, pour lequel la conclusion du théorème de Dussau reste néanmoins valable.

**Mots-clefs** : Conway, Dussau, points cocycliques, droites concourantes, point de Nagel

**Summary**
The famous concyclicity theorem stated by John H. Conway is here reconsidered by means of a parametrisation of the associated triangular configuration with arbitrary triplets of real numbers $(\alpha; \beta; \gamma)$. This theorem, thus corresponding to the case $(\alpha; \beta; \gamma) = (1; 1; 1)$, is generalized while demonstrating that there always exist an infinite family of such triplets which keeps unchanged the conclusion. The « anti-Conway » configuration corresponding to the case $(\alpha; \beta; \gamma) = (-1; -1; -1)$ is also investigated : Xavier Dussau's theorem of concurrent lines is redemonstrated and completed by another concyclicity theorem. It is also proved that there exist in general a unique triplet $(\alpha; \beta; \gamma) \neq (-1; -1; -1)$ which is a function of the sides of the considered triangle and which keeps unchanged the conclusion of Dussau's theorem.

## 1. Introduction

Le mathématicien anglais John H. Conway (1937-2020) est récemment décédé, victime parmi tant d'autres de la pandémie Covid-19. Cet article participe d'un hommage à ce mathématicien génial et prolifique, qu'un groupe encore informel créé en 2020 au sein des IREM de La Réunion et de Mayotte a choisi d'entreprendre.

Nous nous attacherons ici à considérer ce qui constitue indubitablement le plus simple des nombreux « théorèmes de Conway », celui relatif au cercle éponyme : c'est le plus accessible à tous les niveaux d'enseignement, peut-être le plus simple mais pas le moins esthétique. Et nous le ferons en examinant les prolongements qui peuvent en être donnés – lesquels ne sont quant à eux pas aussi élémentaires, loin s'en faut.

La première partie présente ce théorème et ses démonstrations, ainsi que d'autres propriétés de la configuration associée. La seconde partie énonce puis démontre une généralisation du théorème de Conway. La troisième partie montre ensuite comment, en « inversant » les hypothèses faites par Conway, peuvent se déduire deux autres théorèmes tout aussi esthétiques que celui qu'il a énoncé. Une quatrième partie, enfin, s'attache à examiner si, là encore, on peut préserver, tout en généralisant ses hypothèses, la conclusion du premier des théorèmes ainsi obtenu. Un nouveau théorème sera là encore énoncé, puis démontré.

---

[1] Professeur agrégé de mathématiques en CPGE au Lycée Roland Garros du Tampon, La Réunion. L'auteur tient à remercier Gilles Patry, professeur agrégé et IA-IPR de mathématiques de l'Académie de Rennes, pour la réalisation des figures et ses relectures critiques constructives qui ont permis de progresser dans certaines démonstrations. L'auteur remercie également Ivan Riou, professeur agrégé de mathématiques au Centre Universitaire de Mayotte, pour sa relecture attentive de l'article et les détails qu'elle a suscités.

On considère dans toute la suite trois points distincts et non alignés $A, B$ et $C$. On note :
$a = BC, b = AC, c = AB$ et $p$ le demi-périmètre $\frac{1}{2}(a + b + c)$ du triangle $ABC$.
$\Omega$ le centre du cercle inscrit dans $ABC$ et $r$ le rayon du cercle inscrit dans $ABC$.
$(\alpha; \beta; \gamma)$ un triplet de réels quelconques.
$A'$ le point de $(AB)$ tel que $\overrightarrow{AA'} = -\alpha \frac{a}{c} \overrightarrow{AB}$ et $A''$ le point de $(AC)$ tel que $\overrightarrow{AA''} = -\alpha \frac{a}{b} \overrightarrow{AC}$.
$B'$ le point de $(BC)$ tel que $\overrightarrow{BB'} = -\beta \frac{b}{a} \overrightarrow{BC}$ et $B''$ le point de $(BA)$ tel que $\overrightarrow{BB''} = -\beta \frac{b}{c} \overrightarrow{BA}$.
$C'$ le point de $(CA)$ tel que $\overrightarrow{CC'} = -\gamma \frac{c}{b} \overrightarrow{CA}$ et $C''$ le point de $(CB)$ tel que $\overrightarrow{CC''} = -\gamma \frac{c}{a} \overrightarrow{CB}$.

## 2. Le cas $(\alpha; \beta; \gamma) = (1; 1; 1)$ : théorème de Conway et propriétés associées

### 2.1. Le théorème de Conway

Le théorème de Conway peut s'énoncer comme suit avec les notations précédentes :

**Théorème 1** (Conway)
Si $(\alpha; \beta; \gamma) = (1; 1; 1)$, alors les points $A', A'', B', B'', C'$ et $C''$ sont cocycliques ; ils appartiennent à un cercle $(\Gamma)$ dont le centre est le point $\Omega$.

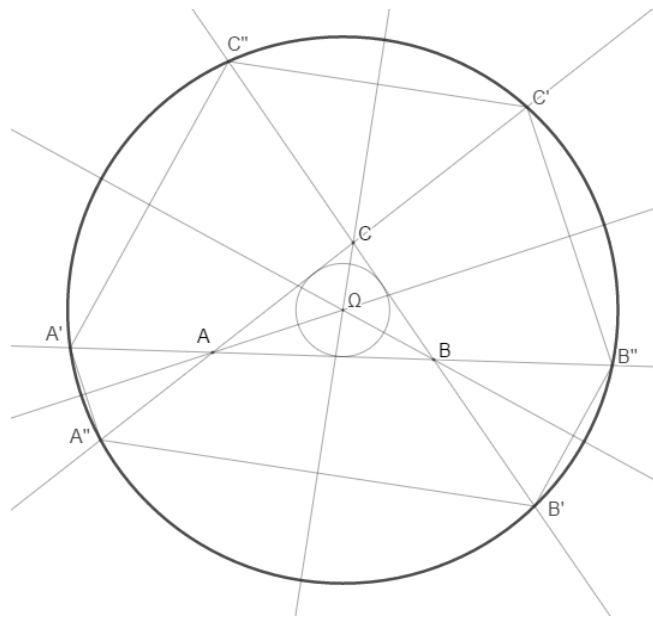

Figure 1

Ce théorème n'a fait l'objet d'aucune publication par Conway lui-même. Il l'a présenté oralement en 1998, probablement pour la première fois, dans un *Mathcamp* (équivalent anglophone de MATH.en.JEANS), sans doute sous la forme d'un défi (celui d'une démonstration) lancé à ses auditeurs[2]. Nous présentons ici deux démonstrations parmi d'autres.

---

[2] Voir par exemple https://blog.tanyakhovanova.com/2010/06/conways-circle/ et http://www.cardcolm.org/JHC.html#Circle

## 2.2 Démonstrations du théorème de Conway

### 2.2.1. Une démonstration purement métrique

Comme $AA' = AA''$, le triangle $AA'A''$ est isocèle en $A$. Donc la bissectrice intérieure du triangle $ABC$ qui est issue de $A$ est la médiatrice de $[A'A'']$. Comme $\Omega$ appartient à cette bissectrice par définition, on a $\Omega A' = \Omega A''$.

Par ailleurs, $BC'' = BC + CC'' = a + c = BA + AA' = BA'$. Donc le triangle $BA'C''$ est isocèle en $B$. Ce qui implique de manière analogue à ce qui précède que $\Omega C'' = \Omega A'$.

Ainsi : $\Omega A' = \Omega A'' = \Omega C''$. Par une méthode analogue, on obtient $\Omega C'' = \Omega C' = \Omega B'' = \Omega B'$. Le théorème de Conway en résulte.

### 2.2.2 Une démonstration angulaire

D'après les hypothèses, $AA'A''$, $BB'B''$ et $CC'C''$ sont des triangles isocèles dont les angles de base sont égaux respectivement à ceux des triangles isocèles $AB''C'$, $BA'C''$ et $CA''B'$. La somme de tous ces angles de base de triangles isocèles (il y en a donc douze) est égale à $4 \times (\widehat{C''C'C} + \widehat{CC'B''} + \widehat{B''B'C''})$. Or, cette même somme n'est autre que celle des six angles au sommet de l'hexagone de Conway ; elle est donc égale à $4 \times 180°$. On en déduit que :

$$\widehat{C''C'C} + \widehat{CC'B''} + \widehat{B''B'C''} = \widehat{C''C'B''} + \widehat{C''B'B''} = 180°$$

Ce qui, comme $[B''C'']$ est une diagonale du quadrilatère convexe $CB''C'C''$, équivaut plus précisément en termes d'angles orientés en radians à $\left(\overrightarrow{C'B''};\overrightarrow{C'C''}\right) = \pi + \left(\overrightarrow{B'B''};\overrightarrow{B'C''}\right)$ $[2\pi]$. La cocyclicité des points $B', B'', C'$ et $C''$ en résulte. Le centre du cercle $(\Gamma)$ qui les contient est l'intersection des médiatrices de $[B'B'']$ et de $[C'C'']$, qui est aussi l'intersection des bissectrices intérieures au triangle $ABC$ respectivement issues de $B$ et de $C$ : c'est donc le point $\Omega$.

On montre par la même méthode que $A', A'', B'$ et $B''$ sont cocycliques ; par un argument analogue, le centre du cercle qui les contient est nécessairement $\Omega$ là encore.

Le théorème de Conway en résulte.

## 2.3. Autres propriétés associées à la configuration de Conway

### 2.3.1. Rayon du cercle de Conway

Notons $R$ le rayon du cercle $(\Gamma)$, $r$ celui du cercle inscrit dans $ABC$ et $U$ le projeté orthogonal de $\Omega$ sur $(AB)$. Alors :

$$R = \sqrt{r^2 + p^2} \text{ , avec la relation classique } r^2 = \frac{(p-a)(p-b)(p-c)}{p}$$

Rappelons que l'expression de $r^2$ découle des deux expressions de l'aire $\mathcal{A}$ de $ABC$, avec la formule de Héron $\mathcal{A} = \sqrt{p(p-a)(p-b)(p-c)}$ et avec la relation $\mathcal{A} = rp$.

Pour la détermination de $R$, on remarque que $A'B'' = A'A + AB + BB'' = a + b + c = 2p$ et que $U$ est le milieu de $[A'B'']$ puisque $\Omega A'B''$ est isocèle en $\Omega$. Le résultat annoncé découle alors immédiatement du théorème de Pythagore appliqué dans $\Omega UA'$ :

$$R^2 = \Omega A'^2 = \Omega U^2 + UA'^2 = r^2 + \left(\frac{1}{2}A'B''\right)^2 = r^2 + p^2$$

### 2.3.2. Propriétés de l'hexagone de Conway

Appelons $A'A''B'B''C'C''$ l'hexagone de Conway déterminé par $(\alpha; \beta; \gamma) = (1; 1; 1)$, qui est clairement convexe. Une propriété en est que ses côtés opposés sont parallèles : ce parallélisme découle immédiatement de trois utilisations du théorème de Thalès (réciproque).

Une autre propriété immédiate de l'hexagone de Conway est que ses trois diagonales $[A'B'']$, $[A''C']$ et $[B'C'']$ sont de même longueur $2p$. De plus, on observe au moyen des réflexions d'axes respectifs $(\Omega A)$, $(\Omega B)$ et $(\Omega C)$ que les six autres diagonales de l'hexagone se regroupent par paires de diagonales de même longueur ($A'C' = A''B''$, $A'B' = B''C''$ et $B'C' = A''C''$).

## 3. Généralisation du théorème de Conway

Nous avons d'emblée généralisé les hypothèses qu'avait faites Conway. Il est naturel dans ces conditions de se demander si la conclusion qu'il énonce dans le cas $(\alpha; \beta; \gamma) = (1; 1; 1)$ demeure valide pour d'autres choix de ce triplet de réels. Nous allons ici établir que la réponse est positive et qu'il y a même une infinité de solutions, que nous allons préciser :

### Théorème 2
On note $\mathcal{T} = \{(\alpha\,;\,1 + (\alpha-1)a/b\,;\,1 + (\alpha-1)a/c) \mid \alpha \in \mathbb{R}\}$.
Les points $A', A'', B', B'', C'$ et $C''$ appartiennent à un même cercle $(\Gamma)$ de centre $\Omega$ si, et seulement si

(i) $ABC$ est scalène et $(\alpha; \beta; \gamma) \in \mathcal{T}$

ou

(ii) $ABC$ est isocèle et $(\alpha; \beta; \gamma) \in \mathcal{T} \cup \left\{\left(0; 0; -\dfrac{a}{c}\right)\right\}$ s'il est isocèle en $C$

$(\alpha; \beta; \gamma) \in \mathcal{T} \cup \left\{\left(0; -\dfrac{c}{b}; 0\right)\right\}$ s'il est isocèle en $B$

$(\alpha; \beta; \gamma) \in \mathcal{T} \cup \left\{\left(-\dfrac{b}{a}; 0; 0\right)\right\}$ s'il est isocèle en $A$

De plus, lorsque cette cocyclicité est réalisée, $U, W, A'$ et $A''$ sont cocycliques, $U, V, B'$ et $B''$ sont cocycliques et $V, W, C'$ et $C''$ sont cocycliques, où les points $U, V$ et $W$ désignent les projetés orthogonaux de $\Omega$ sur les côtés $[AB], [BC]$ et $[CA]$ respectivement.
Le rayon de $(\Gamma)$ est dans tous les cas $R = \sqrt{(p + (\alpha-1)a)^2 + r^2}$.
En particulier, $(\Gamma)$ est le cercle inscrit dans $ABC$ si, et seulement si $\alpha = 1 - p/a$.

Le théorème de Conway apparaît du point de vue du théorème 2 comme un cas très particulier : celui où $\alpha = 1$ dans la classe de solutions formant l'ensemble $\mathcal{T}$.

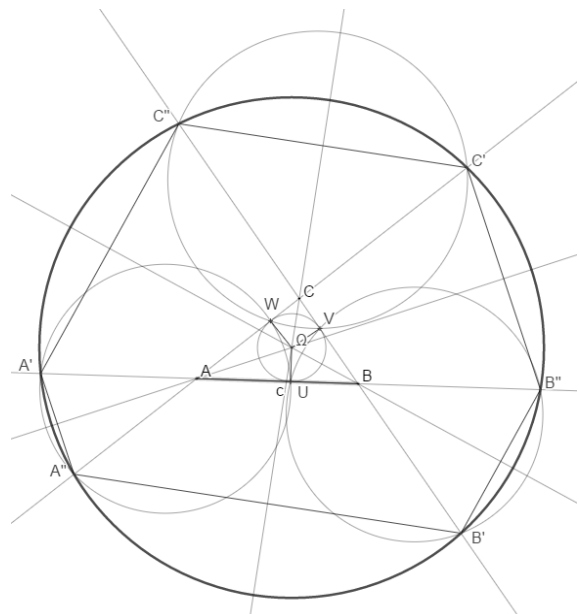

Figure 2

La figure 2 est une illustration de ce théorème général, obtenue avec $\alpha = 3/2$.

On remarquera que cette généralisation prend des libertés avec l'esprit des hypothèses de Conway : si les coefficients sont négatifs, on change en effet le sens des vecteurs. Ceci n'empêche pas la conclusion d'être de même nature, comme le montre la figure 3 avec $\alpha = -1/2$.

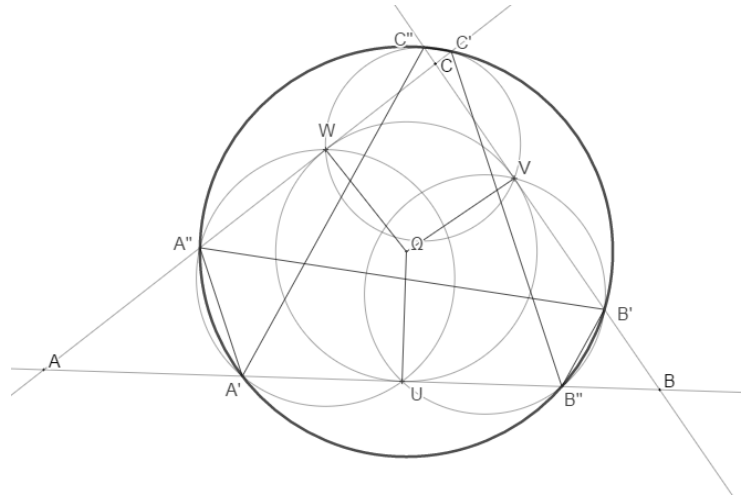

Figure 3

Quant à la dernière assertion du théorème, elle est illustrée par la figure 4, pour un cas où $\alpha$ est voisin de la valeur remarquable $1 - p/a$.

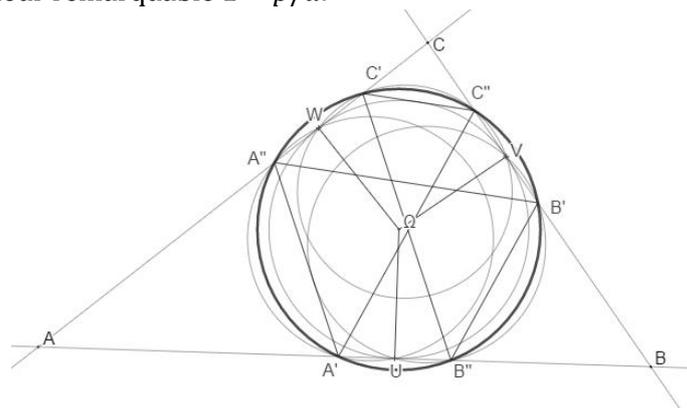

Figure 4

Pour démontrer le théorème 2, nous allons procéder par disjonction de cas, puis à chaque fois « par analyse-synthèse » en supposant d'abord qu'il existe un triplet $(\alpha; \beta; \gamma)$ tel que les points $A', A'', B', B'', C'$ et $C''$ appartiennent à un même cercle $(\Gamma)$ de centre $\Omega$, en trouvant sa forme nécessaire, puis en étudiant réciproquement si cette potentielle solution en est bien une.

Une propriété très importante qui va être plusieurs fois utilisée ici est ce que nous appellerons le « critère de Feuerbach », qui s'énonce comme suit : $D, E, M$ et $N$ étant quatre points distincts du plan tels que les droites $(DE)$ et $(MN)$ ont un point d'intersection $P$, ils sont cocycliques si, et seulement si $\overrightarrow{PD} \cdot \overrightarrow{PE} = \overrightarrow{PM} \cdot \overrightarrow{PN}$.

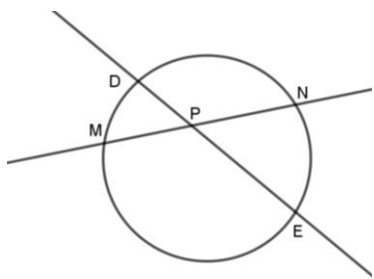 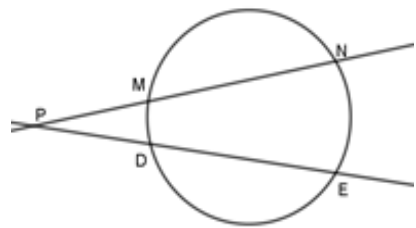

<u>Cas 1</u> : $A' = A''$, $B' = B''$ et $C' = C''$.

Si le cercle $(\Gamma)$ existe, c'est le cercle circonscrit à $ABC$ car on a $A = A' = A''$, $B = B' = B''$ et $C = C' = C''$ ; donc aussi $(\alpha; \beta; \gamma) = (0; 0; 0)$. Mais puisque le centre de $(\Gamma)$ est $\Omega$, il est nécessaire alors que $ABC$ soit équilatéral.

Réciproquement, le triplet $(0; 0; 0)$ est clairement solution sous cette dernière condition, et il est clair qu'il fait dans ce cas partie de l'ensemble $\mathcal{T}$ (on l'obtient avec $a = b = c$ et $\alpha = 0$).

<u>Cas 2</u> : $A' = A''$, $B' = B''$ et $C' \neq C''$. Ce cas implique que $A = A' = A''$ et $B = B' = B''$ ; et $(\alpha; \beta; \gamma) = (0; 0; \gamma)$ avec $\gamma \neq 0$. Distinguons alors à nouveau plusieurs cas.

<u>Cas (2a)</u> : $C', C'', A$ et $B$ sont distincts deux à deux.

On a alors $(AC')$ et $(BC'')$ sécantes en $C$. Si les quatre points sont cocycliques, alors on a l'égalité $\overrightarrow{CA}.\overrightarrow{CC'} = \overrightarrow{CB}.\overrightarrow{CC''}$ d'après le critère de Feuerbach. Or :

$$\overrightarrow{CA}.\overrightarrow{CC'} = \overrightarrow{CA}.\left(-\gamma \frac{c}{b}\overrightarrow{CA}\right) = -\gamma \frac{c}{b}\overrightarrow{CA}^2 = -\gamma bc$$

$$\overrightarrow{CB}.\overrightarrow{CC''} = \overrightarrow{CB}.\left(-\gamma \frac{c}{a}\overrightarrow{CB}\right) = -\gamma \frac{c}{a}\overrightarrow{CB}^2 = -\gamma ac$$

Donc $a = b$. C'est-à-dire que le triangle $ABC$ doit nécessairement être isocèle en $C$. Remarquons qu'on a alors $a > c/2$ d'après l'inégalité triangulaire. Supposons donc maintenant que $a = b$. La puissance de $C$ par rapport au cercle $(\Gamma)$ est $\overrightarrow{CA}.\overrightarrow{CC'} = -\gamma ac = \Omega C^2 - R^2$, avec ici :

$$R^2 = \Omega A^2 = \Omega U^2 + AU^2 = r^2 + \frac{c^2}{4}$$

Par ailleurs, on obtient ici en considérant le triangle $CUA$ que :

$$(\Omega C + r)^2 = CU^2 = CA^2 - AU^2 = a^2 - \frac{c^2}{4}$$

Par conséquent :

$$\Omega C = \frac{1}{2}\sqrt{4a^2 - c^2} - r$$

On a ainsi nécessairement :

$$-\gamma ac = \left(\frac{1}{2}\sqrt{4a^2 - c^2} - r\right)^2 - r^2 - \frac{c^2}{4}$$

Et il n'existe donc qu'une seule solution potentielle :

$$\gamma = \frac{1}{ac}\left(r^2 + \frac{c^2}{4} - \left(\frac{1}{2}\sqrt{4a^2 - c^2} - r\right)^2\right)$$

Cherchons à l'exprimer plus simplement, en tenant compte de :

$$r = \sqrt{\frac{(p-a)^2(p-c)}{p}} = \sqrt{\frac{\left(\frac{c}{2}\right)^2\left(\frac{2a-c}{2}\right)}{\frac{2a+c}{2}}} = \frac{c}{2}\sqrt{\frac{2a-c}{2a+c}}$$

On obtient :

$$\gamma = \frac{1}{ac}\left(\frac{c^2}{4}\frac{2a-c}{2a+c} + \frac{c^2}{4} - \left(\frac{1}{2}\sqrt{4a^2 - c^2} - \frac{c}{2}\sqrt{\frac{2a-c}{2a+c}}\right)^2\right) = \frac{1}{ac}\left(\frac{ac^2}{2a+c} - \left(\frac{a\sqrt{2a-c}}{\sqrt{2a+c}}\right)^2\right)$$

$$= \frac{1}{ac}\left(\frac{ac^2}{2a+c} - \frac{a^2(2a-c)}{2a+c}\right) = \frac{1}{ac}\frac{ac^2 - 2a^3 + a^2c}{2a+c} = \frac{1}{c}\frac{(c^2 - a^2) + a(c-a)}{2a+c}$$

$$= \frac{1}{c} \frac{(c-a)(2a+c)}{2a+c} = 1 - \frac{a}{c}$$

Avec par conséquent aussi $a \neq c$ du fait de $\gamma \neq 0$, ce qui signifie que $ABC$ n'est pas équilatéral.

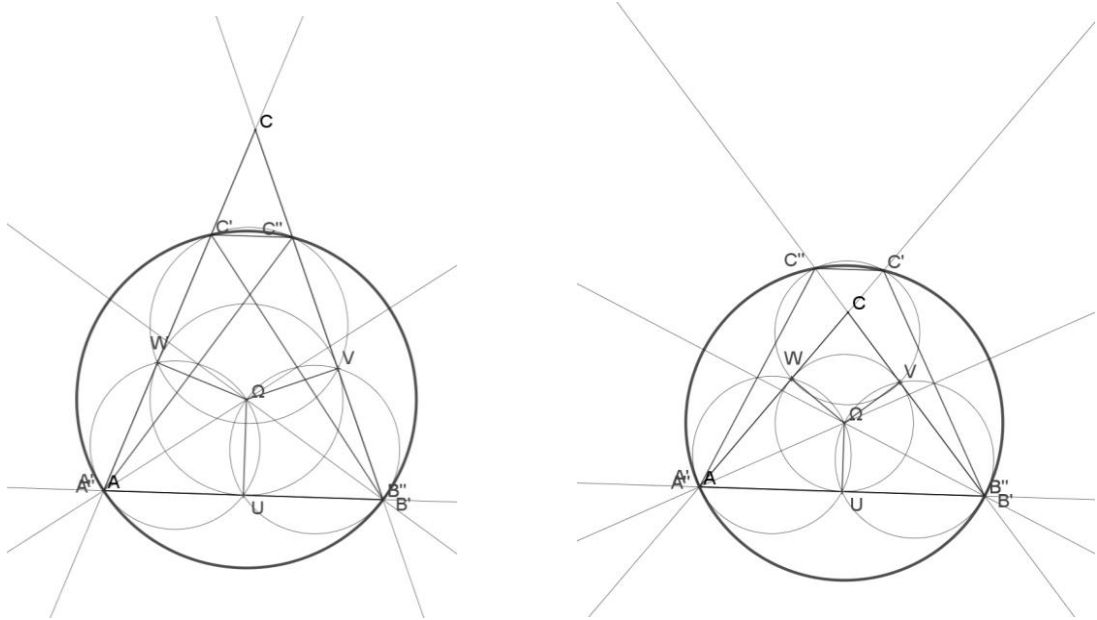

Figures 5 et 5' correspondant au cas (2a), avec $c < a$ ou avec $c > a$

Réciproquement, supposons que $ABC$ est isocèle en $C$ mais non équilatéral et que $(\alpha; \beta; \gamma) = (0; 0; 1 - a/c)$. Alors les points $A, B, C'$ et $C''$ sont cocycliques puisqu'on a $(AC')$ et $(BC'')$ sécantes en $C$ et $\overrightarrow{CA}.\overrightarrow{CC'} = \overrightarrow{CB}.\overrightarrow{CC''} = a(c-a)$. Néanmoins, cela ne montre pas que $\Omega$ est le centre du cercle qui contient ces quatre points. La preuve qui suit va se révéler au cœur de la démonstration aussi dans les cas qui suivent, d'où la mise en évidence d'abord de formules générales, puis leur déclinaison ici dans l'hypothèse $a = b$. Commençons par un peu de trigonométrie. Notons $\theta$ une mesure de l'angle géométrique aigu $\widehat{ACB}/2$. On a :

$$\cos^2(\theta) = \frac{1}{2}(1 + \cos(2\theta)) = \frac{1}{2}\left(1 + \frac{a^2 + b^2 - c^2}{2ab}\right) = \frac{1}{4}\frac{(a+b)^2 - c^2}{ab} = \frac{p(p-c)}{ab}$$

$$\sin^2(\theta) = \frac{1}{2}(1 - \cos(2\theta)) = \frac{1}{2}\left(1 - \frac{a^2 + b^2 - c^2}{2ab}\right) = \frac{1}{4}\frac{c^2 - (a-b)^2}{ab} = \frac{(p-a)(p-b)}{ab}$$

On en déduit en considérant le triangle $\Omega CV$ :

$$\Omega C^2 = \frac{r^2}{\sin^2(\theta)} = \frac{(p-a)(p-b)(p-c)}{p} \frac{ab}{(p-a)(p-b)} = \frac{ab(p-c)}{p}$$

On obtient alors aussi, en considérant cette fois le triangle $C\Omega C'$ (avec pour la deuxième égalité l'englobement en une formule de deux cas à distinguer *a priori*, selon $a \leq c$ ou $a > c$) :

$$\Omega C'^2 = CC'^2 + \Omega C^2 - 2CC' \times \Omega C \times \cos(\widehat{\Omega CC'})$$

$$= \gamma^2 c^2 + \frac{ab(p-c)}{p} + 2\gamma c \sqrt{ab} \frac{\sqrt{p-c}}{\sqrt{p}} \frac{\sqrt{p}\sqrt{p-c}}{\sqrt{ab}} = \gamma^2 c^2 + (p-c)\left(\frac{ab}{p} + 2\gamma c\right)$$

Et donc en particulier ici :

$$\Omega C'^2 = (c-a)^2 + \left(\frac{2a+c}{2} - c\right)\left(\frac{2a^2}{2a+c} + 2\left(1 - \frac{a}{c}\right)c\right) = (c-a)^2 + (2a-c)\left(\frac{a^2}{2a+c} + c - a\right)$$

$$= \frac{(c-a)^2(2a+c) + (2a-c)(c^2 - a^2 + ac)}{2a + c} = \frac{ac^2}{2a+c}$$

Or, on a par ailleurs du fait de l'isocélie en $C$ :

$$\Omega A^2 = r^2 + \frac{c^2}{4} = \frac{(p-a)^2(p-c)}{p} + \frac{c^2}{4} = \frac{\left(\frac{2a+c}{2} - a\right)^2 \left(\frac{2a+c}{2} - c\right)}{\frac{2a+c}{2}} + \frac{c^2}{4} = \frac{ac^2}{2a+c}$$

On a ainsi toujours $\Omega A = \Omega C'$ si $ABC$ est isocèle en $C$ et $(\alpha; \beta; \gamma) = (0; 0; 1 - a/c)$. Ce triplet fournit donc l'unique solution dans la configuration étudiée. On remarque de plus comme au cas 1 que ce triplet appartient à $\mathcal{T}$ (on l'obtient avec $a = b$ et $\alpha = 0$).

<u>Cas (2b)</u> : $C''$, $A$ et $B$ sont distincts, et $C' = A$.

Dans ce cas, on a $\gamma = -b/c$ et donc $CC'' = b$ par définition de $C''$, avec $C''$ dans l'intérieur de $[BC]$ si $b < a$, et sur la demi-droite d'origine $B$ incluse dans $(BC)$ qui ne contient pas $C$ si $a < b$. De plus, si le cercle $(\Gamma)$ existe, alors c'est le cercle circonscrit à $C''AB$. Si l'on suppose que $\Omega$ est le centre de $(\Gamma)$, on a donc $\Omega A = \Omega B = \Omega C''$. Mais puisque $\Omega B \leq \Omega C'' + C''B$, ceci implique que $C''B = 0$, c'est-à-dire $a = b$ et donc $ABC$ isocèle en $C$. On a dans cette situation seulement deux points dans la famille des six qui font l'objet du problème : ce sont $A$ et $B$.

Il est évident réciproquement que si $ABC$ est isocèle en $C$, $\alpha = \beta = 0$ et $\gamma = -b/c = -a/c$, alors on a $C' = A$ et $C'' = B$, de sorte que le cercle centré en $\Omega$ et de rayon $\Omega A = \sqrt{(c/2)^2 + r^2}$ convient. Contrairement aux cas 1 et (2a), on obtient toutefois ici une solution qui n'est pas dans $\mathcal{T}$ : ce type de solution (à permutation près des rôles de $A$, $B$ et $C$) n'existe qu'en cas d'isocélie.

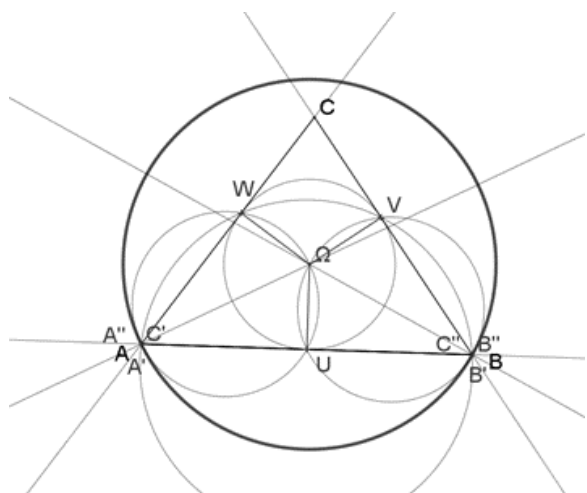

Figure 6 illustrant le cas (2b)

<u>Cas 3</u> : $A' = A''$, $B' \neq B''$ et $C' \neq C''$.

Dans ce cas, $A' = A'' = A$ et on a $\alpha = 0$, $\beta \neq 0$ et $\gamma \neq 0$. Remarquons que si les cinq points $A, B', B'', C'$ et $C''$ sont cocycliques, alors ils appartiennent nécessairement à un cercle $(\Gamma)$ de centre $\Omega$. En effet, les médiatrices respectives de $[B'B'']$ et $[C'C'']$ s'intersectent au centre du cercle circonscrit $(\Gamma)$ au quadrilatère $B''B'C'C''$ ; comme ce sont aussi les bissectrices du triangle $ABC$ respectivement issues de $B$ et de $C$, ce centre est donc nécessairement $\Omega$. Les droites $(AB')$ et $(B''C')$ se coupent en $B$ et les droites $(AC')$ et $(B''C'')$ se coupent en $C$. D'après le critère de Feuerbach, on a donc nécessairement $\begin{cases} \vec{BA}.\vec{BB''} = \vec{BB'}.\vec{BC''} \\ \vec{CA}.\vec{CC'} = \vec{CB'}.\vec{CC''} \end{cases}$. Or :

$$\vec{BA}.\vec{BB''} = \vec{BA}.\left(-\beta \frac{b}{c}\vec{BA}\right) = -\frac{\beta b}{c}\vec{BA}^2 = -\frac{\beta b}{c}c^2 = -\beta bc$$

$$\overrightarrow{BB'}.\overrightarrow{BC''} = \left(-\beta\frac{b}{a}\overrightarrow{BC}\right).\left(\frac{a+\gamma c}{a}\overrightarrow{BC}\right) = -\frac{\beta b(a+\gamma c)}{a^2}\overrightarrow{BC}^2 = -\frac{\beta b(a+\gamma c)}{a^2}a^2 = -\beta b(a+\gamma c)$$

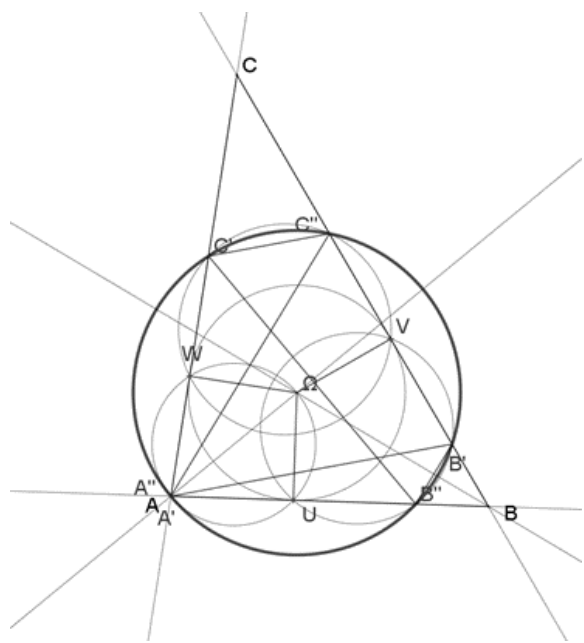

Figure 7 illustrant le cas 3

$$\overrightarrow{CA}.\overrightarrow{CC'} = \overrightarrow{CA}.\left(-\gamma\frac{c}{b}\overrightarrow{CA}\right) = -\gamma\frac{c}{b}\overrightarrow{CA}^2 = -\gamma bc$$

$$\overrightarrow{CB'}.\overrightarrow{CC''} = \left(\frac{a+\beta b}{a}\overrightarrow{CB}\right).\left(-\gamma\frac{c}{a}\overrightarrow{CB}\right) = -\frac{\gamma c(a+\beta b)}{a^2}\overrightarrow{CB}^2 = -\frac{\gamma c(a+\beta b)}{a^2}a^2 = -\gamma c(a+\beta b)$$

Par conséquent : $\begin{cases} \beta bc = \beta b(a+\gamma c) \\ \gamma bc = \gamma c(a+\beta b) \end{cases}$. Ce qui équivaut à $\begin{cases} c = a+\gamma c \\ b = a+\beta b \end{cases}$ puisque $\beta \neq 0$ et $\gamma \neq 0$. On en déduit qu'on a nécessairement $\begin{cases} \beta = 1 - a/b \\ \gamma = 1 - a/c \end{cases}$.

Réciproquement, supposons que $(\alpha;\beta;\gamma) = (0\,;1-a/b\,;1-a/c)$. On a encore ici, comme établi au cas 2a) :

$$\Omega C'^2 = \gamma^2 c^2 + (p-c)\left(\frac{ab}{p} + 2\gamma c\right)$$

La valeur de $\gamma$ fournit donc ici :

$$\Omega C'^2 = (c-a)^2 + (p-c)\left(\frac{ab}{p} + 2c\left(1-\frac{a}{c}\right)\right) = (c-a)^2 + ab - \frac{abc}{p} + 2(c-a)(p-c)$$
$$= (c-a)^2 + ab - \frac{2abc}{a+b+c} + (c-a)(a+b-c) = b(c-a) + ab - \frac{2abc}{a+b+c}$$
$$= bc\left(1 - \frac{2a}{a+b+c}\right) = bc\frac{-a+b+c}{a+b+c}$$

Et on montre de même que :

$$\Omega B'^2 = \beta^2 b^2 + (p-b)\left(\frac{ac}{p} + 2\beta b\right) = bc\frac{-a+b+c}{a+b+c} = \Omega B''^2 = \Omega C''^2$$

Par ailleurs, on obtient comme pour le calcul de $\Omega C^2$ effectué au cas 2a) que :

$$\Omega A^2 = \frac{bc(p-a)}{p} = \frac{2bc}{a+b+c}\cdot\frac{-a+b+c}{2} = bc\frac{-a+b+c}{a+b+c}$$

Ceci établit que le triplet $(\alpha; \beta; \gamma) = (0\,;\,1 - a/b\,;\,1 - a/c)$ est bien dans tous les cas solution, et la seule correspondant à la configuration étudiée. On remarque de plus qu'il appartient lui aussi à $\mathcal{T}$ (il est obtenu avec $\alpha = 0$).

<u>Cas 4</u> : $A' \neq A''$, $B' \neq B''$ et $C' \neq C''$.

Ce cas implique que les neuf points considérés sont deux à deux distincts, et qu'aucun des nombres $\alpha, \beta$ et $\gamma$ n'est nul, avec $(A'B'')$ et $(A''C')$ sécantes en $A$, $(B'C'')$ et $(A'B'')$ sécantes en $B$ et $(C'A'')$ et $(B'C'')$ sécantes en $C$. Si on suppose que les points $A', A'', B', B'', C'$ et $C''$ sont cocycliques, alors d'après le critère de Feuerbach : $\overrightarrow{AA'} \cdot \overrightarrow{AB''} = \overrightarrow{AA''} \cdot \overrightarrow{AC'}$. Or :

$$\overrightarrow{AA'} \cdot \overrightarrow{AB''} = \left(-\alpha \frac{a}{c}\overrightarrow{AB}\right) \cdot \left(\left(1 + \beta \frac{b}{c}\right)\overrightarrow{AB}\right) = -\alpha \frac{a}{c^2}(c + \beta b)\overrightarrow{AB}^2 = -\alpha \frac{a}{c^2}(c + \beta b)c^2 = -\alpha a(c + \beta b)$$

$$\overrightarrow{AA''} \cdot \overrightarrow{AC'} = \left(-\alpha \frac{a}{b}\overrightarrow{AC}\right) \cdot \left(\left(1 + \gamma \frac{c}{b}\right)\overrightarrow{AC}\right) = -\alpha \frac{a}{b^2}(b + \gamma c)\overrightarrow{AC}^2 = -\alpha \frac{a}{b^2}(b + \gamma c)b^2 = -\alpha a(b + \gamma c)$$

Comme $\alpha a \neq 0$, on a donc nécessairement $c + \beta b = b + \gamma c$. On montre de même en utilisant les points $B$ et $C$ qu'on a en fait :

$$\begin{cases} b + \gamma c = c + \beta b \\ c + \alpha a = a + \gamma c \\ a + \beta b = b + \alpha a \end{cases}$$

Or :

$$\begin{cases} b + \gamma c = c + \beta b \\ c + \alpha a = a + \gamma c \\ a + \beta b = b + \alpha a \end{cases} \Leftrightarrow \begin{cases} b\beta - c\gamma = b - c \\ -a\alpha + c\gamma = c - a \\ -a\alpha + b\beta = b - a \end{cases} \Leftrightarrow \begin{cases} -a\alpha + b\beta = b - a \\ -a\alpha + c\gamma = c - a \end{cases} \Leftrightarrow \begin{cases} \beta = 1 + (\alpha - 1)a/b \\ \gamma = 1 + (\alpha - 1)a/c \end{cases}$$

Le triplet $(\alpha; \beta; \gamma) = (\alpha\,;\,1 + (\alpha - 1)a/b\,;\,1 + (\alpha - 1)a/c)$, qui est un élément de $\mathcal{T}$, est donc pour chaque $\alpha \neq 0$ fixé arbitrairement l'unique solution possible ici.

Réciproquement, donnons-nous un triplet $(\alpha; \beta; \gamma)$ de la forme :

$$\left(\alpha\,;\,1 + (\alpha - 1)\frac{a}{b}\,;\,1 + (\alpha - 1)\frac{a}{c}\right), \text{ où } \alpha \in \mathbb{R}^*$$

Toujours par le même calcul que celui effectué au cas 2a) :

$$\Omega C'^2 = \gamma^2 c^2 + (p - c)\left(\frac{ab}{p} + 2\gamma c\right) = \gamma^2 c^2 + ab - \frac{abc}{p} + 2\gamma c(p - c)$$

Donc ici :

$$\Omega C'^2 = \left(1 + (\alpha - 1)\frac{a}{c}\right)^2 c^2 + ab - \frac{2abc}{a + b + c} + 2c\left(1 + (\alpha - 1)\frac{a}{c}\right)\frac{a + b - c}{2}$$

$$= (c + (\alpha - 1)a)^2 + ab - \frac{2abc}{a + b + c} + (c + (\alpha - 1)a)(a + b - c)$$

$$= (c + (\alpha - 1)a)((c + (\alpha - 1)a) + a + b - c) + ab\left(1 - \frac{2c}{a + b + c}\right)$$

$$= \frac{(c + (\alpha - 1)a)(\alpha a + b)(a + b + c) + ab(a + b - c)}{2p}$$

Et on a par ailleurs de même :

$$\Omega A'^2 = \alpha^2 a^2 + bc - \frac{abc}{p} + 2\alpha a(p - a) = \alpha^2 a^2 + bc - \frac{2abc}{a + b + c} + \alpha a(-a + b + c)$$

$$= \alpha a(\alpha a - a + b + c) + bc\left(1 - \frac{2a}{a + b + c}\right) = \alpha a((\alpha - 1)a + b + c) + \frac{bc(-a + b + c)}{a + b + c}$$

$$= \frac{\alpha a(a+b+c)\big((\alpha-1)a+b+c\big)+bc(-a+b+c)}{2p}$$

On vérifie alors facilement après développement et simplifications que :

$$2p(\Omega A'^2 - \Omega C'^2) = \alpha a(a+b+c)\big((\alpha-1)a+b+c\big)+bc(-a+b+c)$$
$$-(c+(\alpha-1)a)(\alpha a+b)(a+b+c)-ab(a+b-c)=0$$

On en déduit que $\Omega A' = \Omega C'$. On montre de la même manière qu'on a en fait aussi bien :

$$\Omega A' = \Omega C' = \Omega B' = \Omega A'' = \Omega B'' = \Omega C''$$

Tout triplet de la forme considérée est donc bien solution du problème, et la seule pour un choix déterminé de $\alpha \neq 0$.

À permutation circulaire près, on a ainsi examiné tous les cas possibles. La première partie du théorème annoncé en résulte. Justifions maintenant les autres assertions.

Les trois cocyclicités annoncées dans le théorème 2 se justifient comme suit. Considérons par exemple le cas de $U, W, A'$ et $A''$. La médiatrice de $[UW]$ est identique à celle de $[A'A'']$ : c'est la droite $(\Omega A)$, par définition des points dans le premier cas, et du fait des résultats précédents dans le second cas. Le quadrilatère $A'A''UW$ est donc un trapèze d'axe de symétrie $(\Omega A)$. On a par conséquent l'égalité angulaire $\widehat{A'WU} = \widehat{A''UW}$ : celle-ci assure la cocyclicité des points $U, W, A'$ et $A''$. Un raisonnement analogue pouvant être fait pour les deux autres quadruplets de points concernés par le résultat annoncé.

Concernant le rayon du cercle $(\Gamma)$, il s'obtient comme suit dans le cas (i) :

$$\overrightarrow{B'C''} = \overrightarrow{B'B} + \overrightarrow{BC} + \overrightarrow{CC''} = \beta\frac{b}{a}\overrightarrow{BC} + \overrightarrow{BC} + \gamma\frac{c}{a}\overrightarrow{BC} = \frac{a+\beta b+\gamma c}{a}\overrightarrow{BC}$$

D'où $B'C'' = |a+\beta b+\gamma c|$, puis $B'V = \frac{1}{2}|a+\beta b+\gamma c|$ puisque $V$ est le milieu de $[B'C'']$ ; de sorte que le théorème de Pythagore fournit, $R$ désignant le rayon de $(\Gamma)$ :

$$R^2 = \Omega B^2 = B'V^2 + \Omega V^2 = \left(\frac{a+\beta b+\gamma c}{2}\right)^2 + r^2$$

Or, on a ici :

$$\frac{a+\beta b+\gamma c}{2} = \frac{1}{2}\left(a+\left(1+(\alpha-1)\frac{a}{b}\right)b+\left(1+(\alpha-1)\frac{a}{c}\right)c\right) = \frac{1}{2}(b+c-a+2\alpha a)$$

$$= \frac{1}{2}(2p+2(\alpha-1)a) = p+(\alpha-1)a$$

On en déduit bien la formule $R = \sqrt{(p+(\alpha-1)a)^2+r^2}$. On peut vérifier qu'elle s'applique aussi dans le cas (ii). En effet, on a vu par exemple en cas d'isocélie en $C$ que $R = \sqrt{(c/2)^2+r^2}$. Or, même dans le cas du triplet $(0; 0; -a/c)$ on obtient :

$$\sqrt{(p+(\alpha-1)a)^2+r^2} = \sqrt{(p-a)^2+r^2} = \sqrt{\left(\frac{2a+c}{2}-a\right)^2+r^2} = \sqrt{\left(\frac{c}{2}\right)^2+r^2}$$

De sorte que la formule obtenue dans le cas (i) s'applique bien encore dans tous les cas du (ii).

Enfin, d'après le précédent résultat, $(\Gamma)$ s'identifie au cercle inscrit dans $ABC$ si, et seulement si $p+(\alpha-1)a=0$. Et cette condition équivaut clairement à $\alpha = 1-p/a$.

Le théorème 2 est ainsi pleinement établi.

## 4. Le cas $(\alpha; \beta; \gamma) = (-1; -1; -1)$ : théorème de Dussau et complément

On peut légitimement appeler le triplet $(-1; -1; -1)$ « anti-Conway » puisque le cas du théorème de Conway correspond au triplet $(1; 1; 1)$, et qualifier de la même manière la configuration qui lui est associée. Or, cette configuration précise est aussi intéressante d'un tout autre point de vue, et elle présente une singularité très esthétique. Les hypothèses générales de cet article vont en fait ici se révéler avoir non pas un, mais deux versants géométriques : non seulement celui de la cocyclicité, que nous avons exploré jusqu'à présent ; mais aussi celui de la concourance de certaines droites… les deux versants se rejoignant justement dans une certaine mesure dans la configuration « anti-Conway ».

### 4.1. La configuration « anti-Conway » et le théorème de Dussau

$X, Y$ et $Z$ étant les points de tangence aux côtés respectifs $[AB], [BC]$ et $[CA]$ des trois cercles exinscrits au triangle $ABC$, il est connu[3] que les droites $(AY), (BZ)$ et $(CX)$ sont concourantes en un point $I$, appelé le point de Nagel de $ABC$.

Considérons alors les six points $A', A'', B', B'', C'$ et $C''$ définis par $(\alpha; \beta; \gamma) = (-1; -1; -1)$. Le résultat associé à cette configuration « anti-Conway » est qu'au lieu d'obtenir la cocyclicité de six points, on obtient en général la concourance de trois droites définies par les six points précédents, conformément au théorème suivant très récemment découvert par le Français Xavier Dussau[4], puis partiellement redécouvert indépendamment par l'auteur du présent article[5] :

**Théorème 3** (Dussau)
Si $ABC$ est scalène et si $(\alpha; \beta; \gamma) = (-1; -1; -1)$, alors les droites $(A'C''), (B'A'')$ et $(C'B'')$ sont concourantes en $I$.

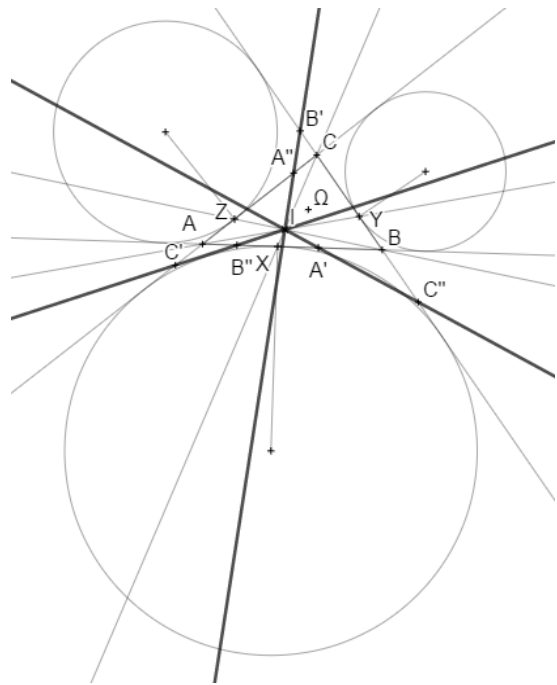

Figure 8

---

[3] Voir par exemple Pamfilos P. : http://users.math.uoc.gr/~pamfilos/eGallery/problems/Nagel.pdf
[4] Dussau X., ''Elementary construction of the Nagel point'', 2020, https://hal.archives-ouvertes.fr/hal-02558108.
[5] L'auteur du présent article, ayant eu la même idée que Dussau d'« inverser » les hypothèses de Conway, a trouvé et démontré par lui-même l'existence de la concourance des droites, sans toutefois observer la coïncidence du point de concours avec le point de Nagel ; et ce fut seulement quelques semaines plus tard.

## 4.2. Démonstration du théorème de Dussau

Les hypothèses permettent d'écrire : $(c-a)\overrightarrow{A'A} + a\overrightarrow{A'B} + 0\overrightarrow{A'C} = \vec{0}$. Ce qui se traduit par l'expression barycentrique : $A' = \text{bar}\{(A; c-a); (B; a); (C; 0)\}$, avec $c \neq a$.

En traduisant de même les cinq autres relations, on obtient aussi, avec $a \neq b$ et $b \neq c$ :

$$B' = \text{bar}\{(A; 0); (B; a-b); (C; b)\} \;;\; C' = \text{bar}\{(A; c); (B; 0); (C; b-c)\}$$
$$A'' = \text{bar}\{(A; b-a); (B; 0); (C; a)\} \;;\; B'' = \text{bar}\{(A; b); (B; c-b); (C; 0)\}$$
$$C'' = \text{bar}\{(A; 0); (B; c); (C; a-c)\}$$

Considérons alors un point $D$ quelconque du plan. Notons $(x; y; z)$ ses coordonnées barycentriques relativement à $A, B$ et $C$ ; c'est-à-dire que $D = \text{bar}\{(A; x); (B; y); (C; z)\}$, avec l'égalité $x + y + z = 1$.

$D$ appartient à la droite $(A'C'')$ si, et seulement si $\begin{vmatrix} c-a & 0 & x \\ a & c & y \\ 0 & a-c & z \end{vmatrix} = 0$, et cette équation est équivalente à $ax + (a-c)y - cz = 0$ puisque $a \neq c$. Par la même méthode, on obtient que $ax - by + (a-b)z = 0$ est une équation barycentrique de $(B'A'')$. L'intersection des droites $(A'C'')$ et $(B'A'')$, si elle existe, a donc des coordonnées barycentriques $(x; y; z)$ solutions du système $\begin{cases} ax + (a-c)y - cz = 0 \\ ax - by + (a-b)z = 0 \end{cases}$. Or, ce système qui est toujours de rang 2 (car $b \neq c$) a une infinité de solutions, donc les deux droites se coupent. Choisissons arbitrairement $z = \frac{a+b-c}{2p}$, qui est non nul du fait de l'inégalité triangulaire. Le système équivaut alors au système de Cramer : $\begin{cases} ax + (a-c)y = c(a+b-c) \\ ax - by = (b-a)(a+b-c) \end{cases}$. Et ses solutions sont données par :

$$\begin{cases} x = \dfrac{\begin{vmatrix} c(a+b-c)/2p & a-c \\ (b-a)(a+b-c)/2p & -b \end{vmatrix}}{\begin{vmatrix} a & a-c \\ a & -b \end{vmatrix}} = \dfrac{-a+b+c}{2p} \\ y = \dfrac{\begin{vmatrix} a & c(a+b-c)/2p \\ a & (b-a)(a+b-c)/2p \end{vmatrix}}{\begin{vmatrix} a & a-c \\ a & -b \end{vmatrix}} = \dfrac{a-b+c}{2p} \end{cases}$$

On en déduit que $(A'C'')$ et $(B'A'')$ se coupent au point $D$ de coordonnées barycentriques $\left(\dfrac{-a+b+c}{2p}; \dfrac{a-b+c}{2p}; \dfrac{a+b-c}{2p}\right)$.

Une équation barycentrique de $(C'B'')$, qui s'obtient par la même méthode que précédemment, est $(b-c)x + by - cz = 0$. Or :

$$(b-c)(-a+b+c) + b(a-b+c) - c(a+b-c) = 0$$

Ceci établit que $D$ appartient aussi à la droite $(C'B'')$.

Les droites $(A'C'')$, $(B'A'')$ et $(C'B'')$ sont donc concourantes en $D$. Reste enfin à observer[6] que le point de Nagel $I$ a précisément les mêmes coordonnées barycentriques que celles trouvées pour $D$, et il est intéressant que Conway ait lui-même connu et commenté le résultat pour $I$ en 1998, sans toutefois observer lui-même l'existence de cette version « inversée » de son théorème[7], qui revient à Dussau et qui consiste à énoncer à la fois l'existence de $D$ et l'égalité $D = I$.

---

[6] Pamfilos P., http://users.math.uoc.gr/~pamfilos/eGallery/problems/Nagel.pdf
[7] http://math.fau.edu/yiu/PSRM2015/yiu/New%20Folder%20(4)/Downloaded%20Papers/barycentrics.pdf

## 4.3. Complément au théorème de Dussau : un autre théorème de cocyclicité

Une propriété complémentaire à celle énoncée par Dussau avec le théorème 3, suggérée par la figure 9, est la suivante :

### Théorème 4
Si $ABC$ est scalène et si $(\alpha;\beta;\gamma)=(-1;-1;-1)$, alors les points $A',A'',B'$ et $C''$ sont cocycliques, les points $B',B'',C'$ et $A''$ sont cocycliques et les points $C',C'',A'$ et $B''$ sont cocycliques.
De plus, le point de Nagel $I$ de $ABC$ a la même puissance par rapport à chacun des cercles circonscrits aux trois quadrilatères correspondants, et cette puissance commune est $8r^2$.

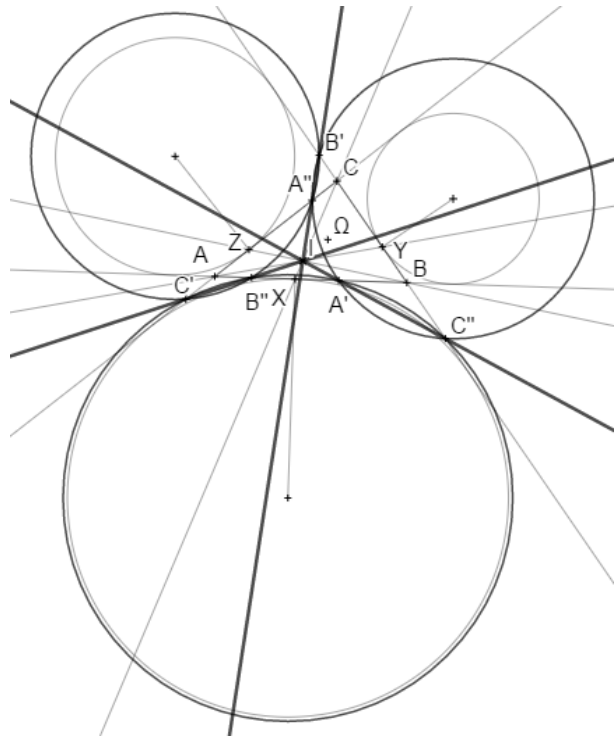

Figure 9

Ici encore, le critère de Feuerbach va se révéler précieux. Nous allons en fait montrer que le calcul du seul produit scalaire $\overrightarrow{IA'}\cdot\overrightarrow{IC''}$ suffit. Qu'il en soit ainsi peut surprendre de prime abord. Mais nous discuterons ensuite ce que cela révèle. Et même si les calculs qui suivent sont rébarbatifs, il est très instructif de voir comment ils révèlent progressivement une parfaite symétrie sous-jacente à la configuration en dépit du caractère scalène du triangle $ABC$, symétrie qu'il nous faudra interpréter adéquatement.

Nous avons vu que :

$$A'=\mathrm{bar}\{(A;c-a);(B;a);(C;0)\} \quad \text{et} \quad I=\mathrm{bar}\{(A;-a+b+c);(B;a-b+c);(C;a+b-c)\}$$

On en déduit les deux relations vectorielles exprimant $\overrightarrow{IA'}$ et $\overrightarrow{IC''}$ dans la base $(\overrightarrow{AB};\overrightarrow{AC})$ des vecteurs du plan $(ABC)$ :

$$\overrightarrow{IA'}=\overrightarrow{AA'}-\overrightarrow{AI}=\left(\frac{a}{c}-\frac{a-b+c}{a+b+c}\right)\overrightarrow{AB}-\frac{a+b-c}{a+b+c}\overrightarrow{AC}=\frac{a+b-c}{c(a+b+c)}\left((a+c)\overrightarrow{AB}-c\overrightarrow{AC}\right)$$

$$=\frac{2(p-c)}{pc}\left((a+c)\overrightarrow{AB}-c\overrightarrow{AC}\right)$$

De même, on déduit de $C'' = \text{bar}\{(A; 0); (B; c); (C; a - c)\}$ que :

$$\overrightarrow{IC''} = \overrightarrow{AC''} - \overrightarrow{AI} = \left(\frac{c}{a} - \frac{a-b+c}{a+b+c}\right)\overrightarrow{AB} + \left(\frac{a-c}{a} - \frac{a+b-c}{a+b+c}\right)\overrightarrow{AC}$$

$$= \frac{-a+b+c}{a(a+b+c)}\left((a+c)\overrightarrow{AB} - c\overrightarrow{AC}\right) = \frac{2(p-a)}{pa}\left((a+c)\overrightarrow{AB} - c\overrightarrow{AC}\right)$$

En tenant compte de $2\overrightarrow{AB}.\overrightarrow{AC} = b^2 + c^2 - a^2$, on obtient donc :

$$\overrightarrow{IA'}.\overrightarrow{IC''} = \frac{4(p-a)(p-c)}{pac}\left((a+c)^2\overrightarrow{AB}^2 - 2c(a+c)\overrightarrow{AB}.\overrightarrow{AC} + c^2\overrightarrow{AC}^2\right)$$

$$= \frac{4(p-a)(p-c)}{pac}\left((c^2(a+c)^2 - c(a+c)(b^2+c^2-a^2) + b^2c^2\right)$$

$$= \frac{4(p-a)(p-c)}{pac}(c(a+c)(a^2-b^2+ac) + b^2c^2)$$

$$= \frac{4(p-a)(p-c)}{p^2ac}ac(a^2+c^2+2ac-b^2) = \frac{4(p-a)(p-c)(a+c-b)(a+c+b)}{p^2}$$

$$= \frac{8(p-a)(p-b)(p-c)}{p}$$

Compte tenu aussi de l'expression de $r$ rappelée au 2.3.1., on en déduit que $\overrightarrow{IA'}.\overrightarrow{IC''} = 8r^2$.

On pourrait alors songer à calculer aussi $\overrightarrow{IB'}.\overrightarrow{IA''}$ pour montrer son égalité au produit scalaire précédent. Mais c'est tout-à-fait superflu. Il est en effet remarquable que le résultat obtenu est une fonction symétrique du triplet $(a; b; c)$, de sorte que l'on peut sans nouveau calcul conclure

$$\overrightarrow{IA'}.\overrightarrow{IC''} = \overrightarrow{IB'}.\overrightarrow{IA''} = \overrightarrow{IC'}.\overrightarrow{IB''} = 8r^2$$

D'une part, les deux premières égalités montrent les trois cocyclicités annoncées en vertu du théorème de Feuerbach. D'autre part, elles expriment la puissance du point $I$ par rapport aux trois cercles circonscrits respectivement aux quadrilatères $A'A''B'C''$, $B'B''C'A''$ et $C'C''A'B''$, la dernière égalité établissant donc que cette puissance est $8r^2$.

Ce résultat interroge quant au statut particulier du point de Nagel du triangle $ABC$. Il semble pertinent de l'interpréter en disant que ce point est un centre de symétrie purement *algébrique* de ce triangle, au sens où c'est un *invariant* sous l'action d'une permutation de l'ensemble $\{A; B; C\}$ ; ce qui s'observe d'ailleurs déjà sur son expression barycentrique. La démonstration précédente révèle progressivement cette symétrie masquée par l'apparente dissymétrie des expressions initiales.

## 5. Extension du théorème de Dussau

Nous avons vu au 3. l'existence d'une infinité de configurations pouvant s'interpréter comme des généralisations de celle qu'avait considéré Conway. Posons ici la question de l'existence analogue d'une généralisation, tout au moins d'une extension du théorème de Dussau.

Adjoignons aux hypothèses initiales de cet article celle que $\alpha, \beta$ et $\gamma$ sont non nuls.

Existe-t-il au moins un triplet $(\alpha; \beta; \gamma) \neq (-1; -1; -1)$ tel que les droites $(A'C'')$, $(B'A'')$ et $(C'B'')$ demeurent concourantes en $I$ ? Nous appellerons ici un tel triplet éventuel une « congruence » de $(-1; -1; -1)$.

La réponse est fournie par le théorème suivant :

**Théorème 5**

Si $ABC$ est scalène, alors :

(i)     Il n'existe aucune congruence de $(-1;-1;-1)$ si $p \in \{2a\,;2b\,;2c\,;\sqrt{2bc}\,;\sqrt{2ca}\,;\sqrt{2ab}\}$.

(ii)    Dans tous les autres cas, il existe une unique congruence de $(-1;-1;-1)$, qui est :

$$(\alpha;\beta;\gamma) = \left(\frac{p^2 - 2bc}{p(p-2a)}\,;\,\frac{p^2 - 2ca}{p(p-2b)}\,;\,\frac{p^2 - 2ab}{p(p-2c)}\right)$$

Aucun des trois nombres de ce triplet congruent n'étant égal à $-1$.

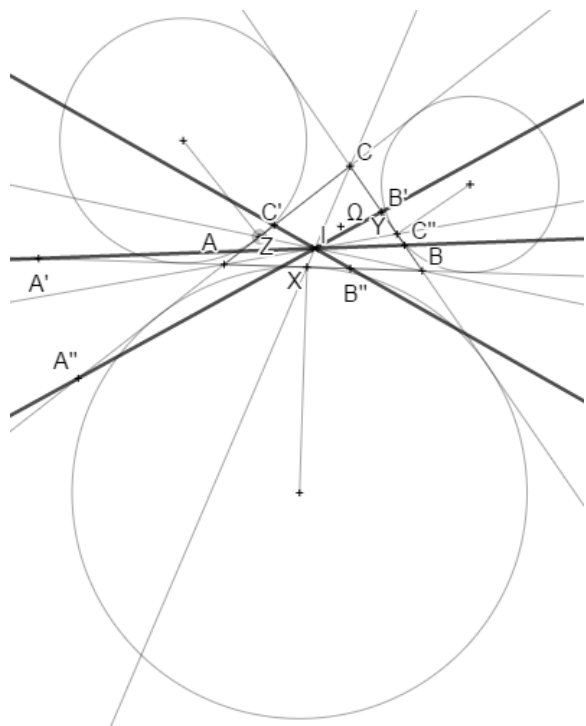

Figure 10

Remarquons au préalable qu'il n'est pas ici exclu d'avoir un (ou deux) des nombres $\alpha$, $\beta$ et $\gamma$ de la solution annoncée au (ii) qui est (ou sont) positif(s). Ce qui est par contre exclu est que ces trois nombres soient simultanément positifs. En effet, on a avec ce triplet $(\alpha;\beta;\gamma)$ particulier :

$$\alpha + \beta + \gamma = \frac{1}{2}(ab + bc + ca - a^2 - b^2 - c^2)$$

$$\leq \frac{1}{2}\left(\frac{1}{2}(a^2 + b^2) + \frac{1}{2}(b^2 + c^2) + \frac{1}{2}(c^2 + a^2) - a^2 - b^2 - c^2\right) = 0$$

avec égalité si, et seulement si $a = b = c$ (résultat qui s'obtient aussi par l'inégalité de Hölder).

Commençons par suivre la même voie que dans la démonstration du théorème 3. On obtient aisément les expressions barycentriques suivantes :

$$A' = \mathrm{bar}\{(A;c+\alpha a); (B;-\alpha a); (C;0)\} \quad ; \quad B' = \mathrm{bar}\{(A;0); (B;a+\beta b); (C;-\beta b)\}$$

$$C' = \mathrm{bar}\{(A;-\gamma c); (B;0); (C;b+\gamma c)\} \quad ; \quad A'' = \mathrm{bar}\{(A;b+\alpha a); (B;0); (C;-\alpha a)\}$$

$$B'' = \mathrm{bar}\{(A;-\beta b); (B;c+\beta b); (C;0)\} \quad ; \quad C'' = \mathrm{bar}\{(A;0); (B;-\gamma c); (C;a+\gamma c)\}$$

Soient $(x;y;z) \in \mathbb{R}^3$ et $D = \mathrm{bar}\{(A;x)\,;(B;y)\,;(C;z)\}$. Il résulte de ce qui précède que :

$$D \in (A'C'') \Leftrightarrow \begin{vmatrix} c+\alpha a & 0 & x \\ -\alpha a & -\gamma c & y \\ 0 & a+\gamma c & z \end{vmatrix} = 0 \Leftrightarrow -\alpha a(a+\gamma c)x - (c+\alpha a)(a+\gamma c)y - \gamma c(c+\alpha a)z = 0$$

De même :
$$D \in (B'A'') \Leftrightarrow -\alpha a(a + \beta b)x - \beta b(b + \alpha a)y - (a + \beta b)(b + \alpha a)z = 0$$
$$D \in (C'B'') \Leftrightarrow -(b + \gamma c)(c + \beta b)x - \beta b(b + \gamma c)y - \gamma c(c + \beta b)z = 0$$

En appliquant ces trois équations barycentriques de droites au point de Nagel que nous avons vu être $I = \text{bar}\{(A; -a + b + c); B(a - b + c); (C; a + b - c)\}$, on obtient que les droites $(A'C''), (B'A'')$ et $(C'B'')$ sont concourantes en $I$ si, et seulement si :

$$(\Sigma) : \begin{cases} a\alpha + c\gamma + p\alpha\gamma = b - p \\ a\alpha + b\beta + p\alpha\beta = c - p \\ b\beta + c\gamma + p\beta\gamma = a - p \end{cases}$$

Système qui est clairement équivalent à :

$$\begin{cases} \alpha(a + p\gamma) = b - p - c\gamma & (E_1) \\ \beta(b + p\gamma) = a - p - c\gamma & (E_2) \\ a\alpha + b\beta + p\alpha\beta = c - p & (E_3) \end{cases}$$

Supposons maintenant que $(\Sigma)$ ait au moins une solution $(\alpha; \beta; \gamma)$ : nous allons ici encore raisonner « par analyse-synthèse ». En fait, on sait déjà qu'il existe une solution, à savoir $(-1; -1; -1)$ ; mais *a priori*, on peut en envisager une qui soit différente. Examinons alors un premier cas : celui où l'on aurait $\gamma = -a/p$. D'après $(E_1)$, on aurait alors $p - b - ac/p = 0$, donc $p^2 - bp - ac = 0$. D'où aussi avec $(E_2)$ :

$$\beta(a - b) = p - a - \frac{ac}{p} = \frac{p^2 - ap - ac}{p} = \frac{(bp + ac) - ap - ac}{p} = \frac{p(b - a)}{p} = b - a$$

Et donc $\beta = -1$ puisque $a \neq b$. En reportant dans $(E_3)$, on en déduirait $\alpha(a - p) = b + c - p$, puis $\alpha(a - p) = 2p - a - p = p - a$, d'où $\alpha = -1$. D'où par conséquent aussi $-a - p\gamma = b - p - c\gamma$ d'après $(E_1)$ ; dont résulterait par un calcul analogue que $\gamma = -1$. En définitive, on aurait donc l'égalité $p = a$, qui est équivalente à $b + c = a$. Ce qui est absurde d'après l'inégalité triangulaire, puisque $ABC$ est non aplati. On montre de même qu'il est impossible d'avoir $\gamma = -b/p$.

On en déduit que s'il existe une solution du système $(\Sigma)$, on a nécessairement :

$$\begin{cases} \alpha = \dfrac{b - p - c\gamma}{a + p\gamma} \\ \beta = \dfrac{a - p - c\gamma}{b + p\gamma} \\ a\dfrac{b - p - c\gamma}{a + p\gamma} + b\dfrac{a - p - c\gamma}{b + p\gamma} - s\dfrac{b - p - c\gamma}{a + p\gamma}\dfrac{a - p - c\gamma}{b + p\gamma} = p - c \quad (E_3') \end{cases}$$

Quelques calculs fournissent alors :

$$(E_3') \Leftrightarrow p(2c - p)\gamma^2 - 2(ab - pc)\gamma + (p^2 - 2ab) = 0$$

On vérifie que $(E_3')$, comme annoncé par le théorème de Dussau, admet $-1$ comme solution :

$$p(2c - p) + 2(ab - pc) + (p^2 - 2ab) = 2pc - p^2 + 2ab - 2pc + p^2 - 2ab = 0$$

L'éventualité de l'existence d'une autre solution dépend alors du degré de cette équation.

Examinons d'abord le cas où $p = 2c$. L'équation $(E_3')$ est alors de degré 1, donc $\gamma = -1$ est sa seule solution. On en déduit aussi, en tenant compte du fait que $b = 3c - a$ ici :

$$\alpha = \frac{b + c - 2c}{a - 2c} = \frac{b - c}{a - 2c} = \frac{(3c - a) - c}{a - 2c} = \frac{2c - a}{a - 2c} = -1$$

De même, on obtient que $\beta = -1$. Il n'existe donc aucune congruence de $(-1; -1; -1)$ dans ce cas, ce dernier triplet étant alors la seule solution de $(\Sigma)$. Le même raisonnement conduit par symétrie à la même conclusion dans les deux cas où $p = 2a$ et $p = 2b$.

Examinons donc maintenant le cas où $p \notin \{2a\,; 2b\,; 2c\,\}$. L'équation $(E_3')$ est alors de degré 2. Et comme $-1$ est l'une de ses solutions, l'opposé du produit de ses solutions est précisément :

$$\gamma = \frac{p^2 - 2ab}{p(p - 2c)}$$

En reportant cette expression dans $(E_1)$ et $(E_2)$, on en déduit respectivement :

$$\alpha = \frac{p^2 - 2bc}{p(p - 2a)} \quad ; \quad \beta = \frac{p^2 - 2ac}{p(p - 2b)}$$

Les trois nombres ainsi obtenus ne sont toutefois solutions du problème posé que s'ils sont non nuls : ceci nécessite donc aussi la condition $p \notin \{\sqrt{2bc}\,; \sqrt{2ca}\,; \sqrt{2ab}\,\}$.

Sous l'hypothèse que $p \notin \{2a\,; 2b\,; 2c\,; \sqrt{2bc}\,; \sqrt{2ca}\,; \sqrt{2ab}\}$, on vérifie réciproquement que le triplet annoncé est bien solution de $(\Sigma)$. Enfin, on sait déjà que $\gamma \neq -1$. Mais on peut observer qu'on a aussi $\alpha \neq -1$ et $\beta \neq -1$. En effet, on a par exemple les implications successives :

$$\alpha = -1 \Rightarrow p^2 - 2bc = 2pa - p^2 \Rightarrow p^2 - ab - pc = 0 \Rightarrow a^2 + b^2 - 2ab - c^2 = 0 \Rightarrow |a - b| = c$$

Ce qui signifie que $a = b + c$ ou $b = a + c$, identités qui sont exclues d'après l'inégalité triangulaire, puisque $ABC$ n'est pas aplati.

Le théorème 5 est ainsi pleinement établi.

## 6. Conclusion

En examinant des modifications des hypothèses faites par Conway, nous avons obtenu plusieurs résultats remarquables dont l'un est une généralisation du théorème qu'il a énoncé.

S'est ainsi manifestée la remarquable puissance de l'algèbre, dans les démonstrations des théorèmes 2, 3 et 4, et plus encore dans celle du théorème 5. Tandis que dans les théorèmes de Conway et de Dussau l'intuition géométrique suffit à imaginer une configuration qui se révèle avoir les propriétés intéressantes qu'énoncent ces théorèmes, c'est en effet l'algèbre qui révèle ici certaines configurations masquées, et qui fait même entièrement le « travail » dans le théorème 5 : nul ne peut en effet imaginer de manière purement géométrique l'existence et la forme de la solution dite « congruente » qu'il annonce, ni ses conditions d'existence ; seul l'algèbre peut en accoucher, et ce sans le secours de la géométrie, qui n'en est que les points de départ et d'arrivée.

S'il est bien difficile de dégager une unité dans l'œuvre foisonnante de Conway, il ne fait aucun doute que ce type de résultat s'accorde impeccablement, tout au moins de ce point de vue d'un dialogue en profondeur entre algèbre et géométrie, à l'esprit flamboyant qui l'a engendrée. Nous espérons avoir ainsi fait œuvre utile à son hommage.